\documentclass[graybox]{svmult}

\usepackage{mathptmx}       
\usepackage{helvet}         
\usepackage{courier}        
\usepackage{type1cm}        
\usepackage{makeidx}         
\usepackage{graphicx}        
\usepackage{multicol}        
\usepackage[bottom]{footmisc}
\usepackage{amsmath}
\usepackage{tabularx}

\newcommand{\sign}{\mathrm{sign}}

\makeindex             

\begin{document}

\title*{Explicit Periodic Solutions in a Delay Differential Equation}
\author{Anatoli Ivanov and Sergiy Shelyag}
\institute{Anatoli Ivanov \at Pennsylvania State University at Wilkes-Barre, 44 University Drive, Dallas, PA 18612, USA, \email{afi1@psu.edu}
\and Sergiy Shelyag \at Flinders University, Tonsley Innovation District, 1284 South Rd, Tonsley, SA, 5042, Australia, \email{sergiy.shelyag@flinders.edu.au}}
%
%
\maketitle

\abstract{We construct stable periodic solutions for a simple form nonlinear delay differential equation (DDE) with a periodic
coefficient. The equation involves one underlying nonlinearity with the multiplicative periodic coefficient.  
The well-known idea of reduction to interval maps is used in the case under consideration, when both the defining nonlinearity 
and the periodic coefficient are piece-wise constant functions. The stable periodic dynamics persist under a smoothing procedure in 
a small neighborhood of the discontinuity set. This work continues the research in recent paper \cite{IvaShe23} on stable periodic
solutions of differential delay equations with periodic coefficients.
}

\section{Introduction}
\label{Intro}

Scalar delay differential equations of the form 
\begin{equation}\label{DDE0}
    x^{\prime}(t)=-\mu x(t)+f(x(t-\tau))
\end{equation}
are used as mathematical models in a broad variety of applications. Some of the most well-known ones are the Mackey-Glass
physiological models \cite{MacGla77}, Nicholson's blowflies models \cite{GurBlyNis80,PerMalCou78}, and many others 
\cite{GlaMac88,Kua93}. 
In spite of its formal simplicity equation (\ref{DDE0}) can exhibit quite complex dynamical behaviors. Unlike scalar ordinary 
differential equations the DDE (\ref{DDE0}) can have complicated dynamics of solutions including oscillations about equilibrium,
stable and unstable periodicity, and complex/chaotic behaviors \cite{DieSvGSVLWal95,HalSVL93}.

More accurate mathematical models of real world phenomena take into account some periodicity factors such as seasonal changes or
circadian rhythm parameter fluctuations and several others. Such improved models with respect to equation (\ref{DDE0}) lead to
similar type DDEs with periodic coefficients:
\begin{equation}\label{DDE1}
    x^{\prime}(t)=-\mu(t) x(t)+a(t) f(x(t-\tau)),
\end{equation}
where $\mu(t)$ and $a(t)$ are periodic functions with the same period. One of the natural questions that appears in this context is
whether DDE (\ref{DDE1}) has periodic solutions with the same period as the periodic coefficients. This work is an attempt to answer
such question in a particular case. 

Note that the periodicity problem we consider is different from similar ones in many other publications, such as in e.g. \cite{Far17} 
(also see additional references therein). We assume the negative feedback conditions on $f, x\cdot f(x)<0,\;\forall x\ne0,$ implying
that DDE (\ref{DDE1}) admits the trivial solutions $x\equiv0$ (for any functions $\mu$ and $a$). This situation is similar to the 
case of autonomous equation (\ref{DDE0}), when it admits a constant solution (with the negative feedback in many cases). In
\cite{Far17} and elsewhere, when such constant positive equilibrium is perturbed by a multiplicative periodic coefficient, the new
problem does not admit constant solutions (positive equilibrium) any longer. 

A particular case of equation (\ref{DDE1}) is considered when $\mu(t)\equiv0$ (see equation (\ref{DDE}) below). The initial choice of
the nonlinearity $f$ and the periodic coefficient $a(t)$ is as piecewise constant functions. Such  simple form allows for an 
explicit calculation of solutions for all forward times as piecewise affine functions (made up continuously of line segments). 
The problem of existence of slowly oscillating solutions and their stability is reduced to the existence of fixed points of a 
simple interval map and their attractivity. The dynamics of both is shown to persist when the two discontinuous functions are
replaced by close to them continuous functions (or even of $C^\infty$ class).

The results of this paper are a continuation and an expansion of recently obtained similar results by the same authors 
\cite{IvaShe23}. 

\section{Preliminaries}
\label{Prelim}

As in \cite{IvaShe23} we consider the special case of DDE (\ref{DDE1}) when $\mu(t)\equiv0$:
\begin{equation}\label{DDE}
    x'(t)=a(t)f(x(t-1))
    \end{equation}
Here initially the nonlinear function $f$ is given by $f(x)=f_0(x)=-\sign(x)$ and the periodic coefficient $a=a(t,a_1,a_2,p_1,p_2)$ 
is defined as
\begin{equation}
a(t) = a_0(t) = \begin{cases}
a_1,\;\text{if}\; 0\le t < p_1 \\
a_2,\;\text{if}\; p_1\le t < p_1+p_2\\
\text{periodic extension outside}\; [0,p_1+p_2)\;\text{for all}\; t\in\mathbf R,
\end{cases}
\label{DDE_COEF}
\end{equation}
with all the values $a_1, a_2, p_1, p_2$ being positive and $p_1+p_2:=T > 1$.

For arbitrary initial function $\varphi(s)\in C=C([-1,0],\mathbf{R})$ such that $\varphi>0\; \forall s\in[-1,0]$ the corresponding 
unique solution $x=x(t,\varphi)$ exists for all $t>0$ (it is calculated by the consecutive integration, the "step method"). 
The solution is a continuous piece-wise affine function made up of straight line segments. It is differentiable everywhere except a 
countable set of points where two line segments with different slopes match.

\section{Main Results}
\label{Main}

\subsection{Explicit Periodic Solutions}
\label{ExplicitPS}

Given $f=f_0$ and $a=a_0$ as above for arbitrary initial function $\phi\in C$ the corresponding solution $x=x(t,\phi)$ is
constructed for all forward times $t\ge0$. It is a continuous piecewise differentiable function composed of consecutive segments 
of affine (linear) functions.

We assume that periodic solutions have the form as depicted in Fig. 1 (this will be shown to be the case for particular values of the
parameters $a_1,a_2,p_1,p_2$). For arbitrary initial function $\phi\in C$ such that $\phi(s)>0\; \forall s\in[-1,0],$ the corresponding
forward solution $x(t),t\ge0,$ depends only on the value $\phi(0):=h>0$ and does not depend on the other values $\phi(s), s\in[-1,0)$
of the initial interval. There exist consecutive zeros $t_1>0$ and $t_2=t_1+2$ of the solution such that $x(t)>0, t\in[0,t_1)$, 
$x(t)<0, t\in(t_1,t_1+2)$, and $x(t)>0, t\in[t_1+2,p_1]$. The necessary condition for such first two zeros to exist is that $p_1>2$.
We also assume that $p_1<t_1+3<p_1+p_2$ and that the solution remains positive on the interval $[p_1,p_1+p_2]$. For the latter to be 
true one only has to require that $x_3:=x(p_1+p_2)>0$.

\begin{figure*}
    \centering
    \includegraphics[width=0.9\textwidth]{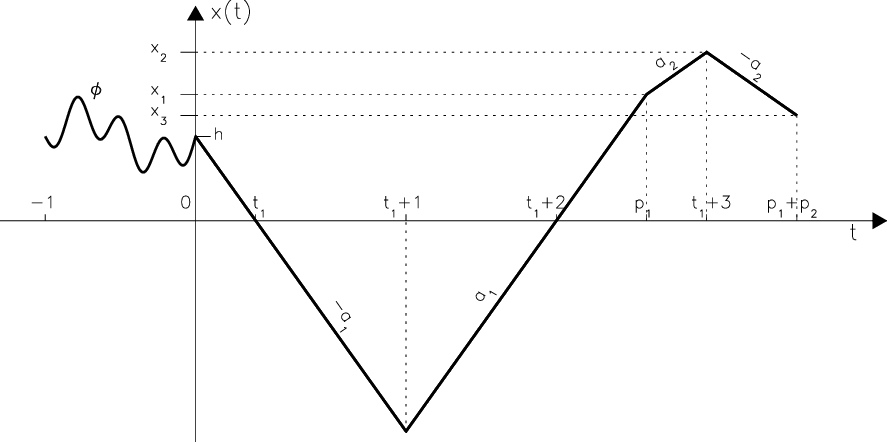}
    \caption{Slowly oscillating piece-wise affine solution}
    \label{fig:enter-label}
\end{figure*}

We next calculate the explicit values of $t_1, x_1=x(p_1), x_2=x(t_1+3), x_3=x(p_1+p_2)$ in terms of the parameters $a_1,a_2, p_1, p_2$ 
and given initial value $h>0$. It is straightforward to find that:
$$
t_1=\frac{h}{a_1},\quad x_1=x(p_1)=-h+a_1p_1-2a_1,\; 
$$
$$
x_2=x(t_1+3)=\left(\frac{a_2}{a_1}-1\right)h+a_1p_1-2a_1+3a_2-a_2p_1,
$$
$$
x_3=x(p_1+p_2)=\left(2\frac{a_2}{a_1}-1\right)h+a_1(p_1-2)+a_2[6-(2p_1+p_2)]:=mh+b
$$
In addition to the previous assumptions the condition $b>0$, i.e.
\begin{equation}\label{b}
    b=a_1(p_1-2)+a_2[6-(2p_1+p_2)]>0,
\end{equation}
guaranties that $x_3>0$, and thus that the solution $x(t), t\in[0,p_1+p_2],$ is of desired shape (as shown in Fig. 1).

A fixed point $h_*>0$ of the affine map $F(x)=mx+b$ gives rise to a slowly oscillating periodic solution $x=x(t,h_*)$ of equation 
(\ref{DDE0}). Moreover, by the linearity of $F$ and the construction, such periodic solution is asymptotically stable if $|m|<1$. 
The latter is equivalent to $0<a_2<a_1$.
The unique fixed point $x=h_*$ is easily found as 
\begin{equation}\label{h}
    h_*=b/(1-m)=\frac{a_1[a_1(p_1-2)+a_2[6-(2p_1+p_2)]]}{2(a_1-a_2)}.
\end{equation}
Therefore, we arrive at the following statement
\begin{theorem}\label{Thm1}
Suppose that the parameters $a_1,a_2,p_1,p_2$ are such that the inequality (\ref{b}) is satisfied and $a_1>a_2$. 
Then DDE (\ref{DDE}) has an asymptotically stable slowly oscillating periodic solution. The periodic solution is generated by the 
initial function $\phi(s)\equiv h_*, s\in[-1,0],$ where $h_*$ is given by (\ref{h}).
\end{theorem}

Note that due to the symmetry property of the nonlinearity $f$ (oddness) and the procedure of construction of the periodic solution
$x_{h_*}$, under the assumption of Theorem \ref{Thm1} there also exists the symmetric to $x_{h_*}$ periodic solution generated by
the initial function $\phi(s)\equiv -h_*$. The two periodic solutions are related by $x_{-h_*}(t)\equiv-x_{h_*}(t)$.

\subsection{Smoothing of Nonlinearities}
\label{Smooth}

In this subsection we demonstrate that the stable periodic solution derived in the previous subsection persists under the standard
smoothing procedure for functions $f$ and $a$. The procedure is a well known one which  has been used in many publications 
(see e.g. \cite{IvaSha91,Pet80} and further references therein). 
It consists in replacing each jump discontinuity by a continuous affine function in a small neighborhood of a discontinuity point. 
We repeat here the exact constructions as in \cite{IvaShe23}.

Consider the functions $f$ and $a$ as defined in subsection \ref{ExplicitPS} : $f(x)=f_0(x)$ and $a(t)=a_0(t)$.
Let $\delta_0>0$ be small, and for every $\delta\in(0,\delta_0]$ introduce the continuous functions
$f_\delta(x)$ and $a_0^\delta(t)$ by:
\begin{equation}\label{f_d}
    f(x)=f_{\delta}(x)=
\begin{cases}
+1\; \text{if}\; x\le -\delta\\
-1\; \text{if}\; x\ge \delta\\
-({1}/{\delta}) x\; \text{if}\; x\in[-\delta, \delta],
\end{cases}
\end{equation}
and
\begin{equation}\label{A0d}
a(t)=a_{\delta}(t)=\begin{cases}
a_2+\frac{a_1-a_2}{2\delta}(t+\delta)\;\text{if}\; t\in[-\delta,\delta]\\
a_1\;\text{if}\; t\in[\delta,p_1-\delta)\\
a_1+\frac{a_2-a_1}{2\delta}[t-(p_1-\delta)]\;\text{if}\; t\in[p_1-\delta,p_1+\delta]\\
a_2\;\text{if}\; t\in[p_1+\delta,p_1+p_2-\delta)\\
a_2+\frac{a_1-a_2}{2\delta}[t-(p_2-\delta)]\;\text{if}\; t\in[p_1+p_2-\delta,p_1+p_2+\delta]\\
\text{periodic extension on}\; \mathbf{R}\; \text{outside interval}\; [0,p_1+p_2).
\end{cases}
\end{equation}

The following statement is an extension of the main Theorem \ref{Thm1} to the case of continuous functions $f_\delta$ and $a_\delta$.

\begin{theorem}\label{Thm2}
Suppose the assumptions of Theorem \ref{Thm1} are satisfied. There exists $\delta_0>0$ such that for every $\delta\in[0,\delta_0]$ 
DDE (\ref{DDE}) with $f=f_\delta$ and $a=a_\delta$ has an asymptotically stable slowly oscillating periodic solution $x_\delta(t).$
The solution $x_\delta(t)$ is close in the uniform metric to the periodic solution $x_{h_*}$ of Theorem \ref{Thm1} for small
$\delta\ge0$, moreover
$\lim_{\delta\to0}\{\sup_{t\in\mathbf{R}}|x_\delta(t)-x_{h_*}(t)|\}=0$.
\end{theorem}

The proof is based on the explicit calculation of the similar value $x_3$ of the solution $x(t,h)$ of the corresponding DDE 
(\ref{DDE}). Its shape is similar to that shown in Fig. 1, except that the corner points (where the derivative $x^\prime$ is not
defined) are replaced by a smooth $C^1$ curve, based on the affine representation of both $f$ and $a$ around their discontinuity
points. The exact calculations repeat those in \cite{IvaShe23}. The final result is a representation of $x_3(h)$ in the form 
$x_3(h)=\Tilde{F}(h)$, where $\Tilde{F}(h)$ is $C^1$-close to the $F$ defined in subsection \ref{ExplicitPS} (see pages 74-75 of 
\cite{IvaShe23} for more details). Therefore, the existence, stability, and closeness follow.

\subsection{Numerical Demonstration }
\label{Numer}

Periodic solutions described in subsections \ref{ExplicitPS} and \ref{Smooth} can be easily verified numerically. It is clear that 
sets of parameters leading to the existence of a stable periodic solution, as described by  Theorem \ref{Thm1}, are not 
exceptional. Any sufficiently small perturbation of a particular quadruple of the parameters will produce an asymptotically stable
periodic solution. They are open in the set of all parameters. The exemplary parametric values for which numerical periodic solutions
have been obtained and confirmed to be of the described form are given in Table~\ref{tab:param_table}.

\begin{table}
    \centering
    \begin{tabular}{|c|c|c|c|c|c|}
\hline    
$~~~~~a_1~~~~~$ & $~~~~~a_2~~~~~$ & $~~~~~p_1~~~~~$ & $~~~~~p_2~~~~~$ & $~~~~~h_*~~~~~$ & $~~~~~T~~~~~$ \\
\hline
1 &	0.25	& 2.5	& 1.5	& 0.25	& 4 \\
2	& 0.5	& 2.5	& 2	& 1/3	& 4.5 \\
2	& 0.25	& 2.5 &	1 &	4/7	 & 3.5 \\
1	& 0.5	& 3	& 1 &	0.5	& 4 \\
2	& 1	& 3	& 1.5 &	0.5 &	4.5 \\
2.5	& 0.5	& 3	& 4 & 0.31 & 	7 \\
3	& 0.5	& 3	& 4.5	& 0.45	& 7.5 \\
5	& 0.5	& 3	& 3	& 1.94 &	6 \\
5	& 1	& 3	& 2	& 1.88	& 5 \\
\hline
    \end{tabular}
    \caption{Parametric values for which stable periodic solutions of equation~(\ref{DDE}) with the coefficient defined in equation~(\ref{DDE_COEF}) have been demonstrated numerically.}
    \label{tab:param_table}
\end{table}

%
%


\vspace{-0.5cm}

\section{Discussion}
\label{Disc}

This paper demonstrates the existence of another type of stable periodic solutions in DDE (\ref{DDE}). They are complimentary to 
those which existence was shown in our recent work \cite{IvaShe23}. There is a possibility of existence of additional types of 
stable periodic solutions in (\ref{DDE}) with piece-wise constant defining functions $f$ and $a$, the question we intend to tackle 
in our next research considerations. It is of interest to show the existence of stable periodic solutions in a general DDE of the form 
(\ref{DDE}), where $f$ and $a$ are smooth functions of arbitrary type. An interesting and seemingly very nontrivial problem is the
existence of (stable) periodic solutions in DDE (\ref{DDE}) with smooth functions, when $f$ is arbitrary and fixed nonlinearity with
the negative feedback property and the periodic coefficient $a$ is parameter dependent, e.g., $a=a(t,\varepsilon)$, such that 
$a(t,0)=a_0>0$ and $a(t+T,\varepsilon)=a(t,\varepsilon)$, for some $T>0$. As the value $\varepsilon=0$ gives periodic solutions in the
well-known autonomous case of equation (\ref{DDE0}), and some larger positive values of $\varepsilon_0$ result in $T$-periodic 
solutions with the same period as $a$, the intermediate values of $\varepsilon\in(0,\varepsilon_0)$ may yield varied dynamics and 
complex transition patterns. This is an open problem worthy further investigation.

\begin{acknowledgement}
The authors thank the mathematical research institute MATRIX in Australia where part of this research was performed. Its final version
resulted from collaborative activities of the authors during the workshop "Delay Differential Equations and Their Applications" 
(https://www.matrix-inst.org.au/events/delay-differential-equations-and-their-applications/) held in December 2023. The authors are 
also grateful for the financial support provided for these research activities by Simons Foundation (USA), Flinders University 
(Australia), and the Pennsylvania State University (USA).
A.I.'s research was also supported in part by the Alexander von Humboldt Stiftung (Germany) during his visit to 
Justus-Liebig-Universit\"{a}t, Giessen, in June-August 2023.
\end{acknowledgement}
%

%
%
%

\end{document}